\documentclass[oneside,notitlepage,10pt]{amsart}
\usepackage{geometry}
\usepackage{amsfonts,amsmath,latexsym,amssymb}
\usepackage{amsmath,amsfonts,amssymb}
\usepackage{graphicx}
\usepackage{epstopdf}
\usepackage[colorlinks=true,urlcolor=magenta,linkcolor=red,citecolor=blue]{hyperref}

\newtheorem{theorem}{Theorem}
\newtheorem{lemma}{Lemma}
\newtheorem{corollary}{Corollary}
\theoremstyle{definition}

\newtheorem{remark}{Remark}

\hypersetup{
pdftitle={New Identities on the Apostol-Bernoulli Polynomials of Higher Order},
pdfauthor={\textcopyright\ Armen Bagdasaryan, Serkan Araci},
}
\input{tcilatex}

\begin{document}
\title[New Identities on the Apostol-Bernoulli Polynomials]{New Identities
on the Apostol-Bernoulli Polynomials of Higher Order Derived From Bernoulli
Basis}
\author{Armen Bagdasaryan}
\address{Department of Mathematics and Statistics, American University of
the Middle East, Kuwait City, 15453 Egaila, Kuwait; and Russian Academy of
Sciences, Institute for Control Sciences, 65 Profsoyuznaya, 117997 Moscow,
Russia; }
\email{bagdasar@member.ams.org}
\author{Serkan Araci}
\address{Department of Economics, Faculty of Economics, Administrative and
Social Sciences, Hasan Kalyoncu University\\
TR-27410 Gaziantep, Turkey}
\email{mtsrkn@hotmail.com}
\author{Mehmet Acikgoz}
\address{Department of Mathematics, Faculty of Arts and Sciences, University
of Gaziantep TR-27310 Gaziantep, Turkey}
\email{acikgoz@gantep.edu.tr}
\author{Yuan He}
\address{Faculty of Science, Kunming University of Science and Technology,
Kunming, Yunnan 650500, People's Republic of China}
\email{hyyhe@aliyun.com}
\subjclass[2000]{Primary 11B68, 11S80}
\date{January 20, 2015}
\keywords{Generating function; Bernoulli polynomials of higher order; Euler
polynomials of higher order; Hermite polynomials; Apostol-Bernoulli
polynomials of higher order; Apostol-Euler polynomials of higher order;
Identities}

\begin{abstract}
In the present paper, we obtain new interesting relations and identities of
the Apostol-Bernoulli polynomials of higher order, which are derived using a
Bernoulli polynomial basis. Finally, by utilizing our method, we also derive
formulas for the convolutions of Bernoulli and Euler polynomials, expressed
via Apostol-Bernoulli polynomials of higher order.
\end{abstract}

\maketitle












\section{Introduction}

For $t\in \mathbb{C}$, the Euler polynomials have the following Taylor
expansion at $t=0$:%
\begin{equation}
\sum_{n=0}^{\infty }E_{n}\left( x\right) \frac{t^{n}}{n!}=e^{tE\left(
x\right) }=\frac{2}{e^{t}+1}e^{xt}\text{, }\left( \left\vert t\right\vert
<\pi \right)  \label{eqn1}
\end{equation}%
with the usual convention about replacing of $\left( E\left( x\right)
\right) ^{n}:=E_{n}\left( x\right) $ (see \cite{Araci6}, \cite{Jolany}, \cite%
{He1}, \cite{Kim6}, \cite{Kim8}, \cite{Luo1}, \cite{Srivastava1}, \cite%
{Srivastava2}).%

There are also explicit formulas for the Euler polynomials, e.g. 
\begin{equation*}
E_{n}(x)=\sum_{k=0}^{n}\binom{n}{k}\frac{E_{k}}{2^{k}}\left( x-\frac{1}{2}%
\right) ^{n-k}
\end{equation*}%
where $E_{k}$ means the Euler numbers. Conversely, the Euler numbers are
expressed with the Euler polynomials through $E_{k}=2^{k}E_{k}(1/2)$. These
numbers can be computed by:%
\begin{equation*}
\left( E+1\right) ^{n}+\left( E-1\right) ^{n}=\left\{ 
\begin{array}{cc}
2 & \text{if }n=0 \\ 
0 & \text{if }n\neq 0%
\end{array}%
\right.
\end{equation*}%
(see \cite{Hansen}, \cite{Norlund}).

For $\left\vert t\right\vert <2\pi $ with $t\in \mathbb{C}$, the Bernoulli
polynomials are defined by%
\begin{equation*}
\sum_{n=0}^{\infty }B_{n}\left( x\right) \frac{t^{n}}{n!}=e^{tB\left(
x\right) }=\frac{t}{e^{t}-1}e^{xt}
\end{equation*}%
where we have used $\left( B\left( x\right) \right) ^{n}:=B_{n}\left(
x\right) $, symbolically. In the case $x=0$, we have $B_{n}\left( 0\right)
:=B_{n}$ that stands for $n$-th Bernoulli number. This number can be
computed via%
\begin{equation*}
\left( B+1\right) ^{n}-B_{n}=\delta _{n,1}
\end{equation*}%
where $\delta _{n,1}$ stands for Kronecker delta (see \cite{Bagdasaryan}, 
\cite{dilcher}, \cite{Kim3}, \cite{Kim6}).

The Euler polynomials of order $k$ are defined by the exponential generating
function as follows: 
\begin{equation}
\left( \frac{2}{e^{t}+1}\right)
^{k}e^{xt}=e^{tE^{(k)}(x)}=\sum_{n=0}^{\infty }E_{n}^{(k)}(x)\frac{t^{n}}{n!}%
\quad (k\in \mathbb{Z}_{+}=\mathbb{N}\cup \left\{ 0\right\} )\text{,}
\label{eqn2}
\end{equation}%
with the usual convention about replacing $(E^{(k)}(x))^{n}$ by $%
E_{n}^{(k)}(x)$. In the special case, $x=0,\,E_{n}^{(k)}(0):=E_{n}^{(k)}$
are called Apostol-Euler numbers of order $k$ (see \cite{Kim8}, \cite{Luo1}).

In the complex plane, Apostol-Euler polynomials $E_{n}\left( x\mid \lambda
\right) $ and Apostol-Bernoulli polynomials $B_{n}\left( x\mid \lambda
\right) $ are given by \cite{Luo1}%
\begin{eqnarray}
\frac{2}{\lambda e^{t}+1}e^{xt} &=&\sum_{n=0}^{\infty }E_{n}\left( x\mid
\lambda \right) \frac{t^{n}}{n!},\text{ }\left( \left\vert t\right\vert
<\log \left( -\lambda \right) \right)  \label{eqn3} \\
\frac{t}{\lambda e^{t}-1}e^{xt} &=&\sum_{n=0}^{\infty }B_{n}\left( x\mid
\lambda \right) \frac{t^{n}}{n!},\text{ }\left( \left\vert t\right\vert
<\log \lambda \right) \text{.}  \label{equat4}
\end{eqnarray}

In \cite{Luo1}, Apostol-Euler polynomials of higher order $E_{n}^{\left(
k\right) }\left( x\mid \lambda \right) $ and Apostol-Bernoulli polynomials
of higher order $B_{n}^{\left( k\right) }\left( x\mid \lambda \right) $ are
given by the following generating functions:%
\begin{eqnarray}
\left( \frac{2}{\lambda e^{t}+1}\right) ^{k}e^{xt} &=&\sum_{n=0}^{\infty
}E_{n}^{\left( k\right) }\left( x\mid \lambda \right) \frac{t^{n}}{n!}%
,\left( \left\vert t\right\vert <\log \left( -\lambda \right) \right)
\label{eqn4} \\
\frac{t^{k}}{\left( \lambda e^{t}-1\right) ^{k}}e^{xt} &=&\sum_{n=0}^{\infty
}B_{n}^{\left( k\right) }\left( x\mid \lambda \right) \frac{t^{n}}{n!}%
,\left( \left\vert t\right\vert <\log \lambda \right) \text{.}
\label{eqaution7}
\end{eqnarray}

In the above expressions, we take the principal value of the logarithm $\log
\lambda $, i.e., $\log \lambda =\log \left\vert \lambda \right\vert +i\arg
\lambda \left( -\pi <\arg \lambda \leq \pi \right) $ when $\lambda \neq 1;$
set $\log 1=0$ when $\lambda =1.$ Additionally, in the special case, $x=0$
or $\lambda=1$ in (\ref{eqn4}) and (\ref{eqaution7}), we have $E_{n}^{\left(
k\right) }\left( 0\mid \lambda \right) :=E_{n}^{\left( k\right) }\left(
\lambda \right) $ and $B_{n}^{\left( k\right) }\left( 0\mid \lambda \right)
:=B_{n}^{\left( k\right) }\left( \lambda \right)$, $E_{n}^{\left( k\right)
}\left( x\mid 1 \right) :=E_{n}^{\left( k\right) }\left( x\right) $ and $%
B_{n}^{\left( k\right) }\left( x\mid 1\right) :=B_{n}^{\left( k\right)
}\left( x \right) $ that stand for Apostol-Euler numbers, Apostol-Bernoulli
numbers, the Euler polynomials of order $k$ and the Bernoulli polynomials of
order $k$.

Apostol-Euler polynomials of higher order and Apostol-Bernoulli polynomials
of higher order can be expressed in terms of their numbers as follows:%
\begin{equation}
E_{n}^{\left( k\right) }\left( x\mid \lambda \right) =\sum_{l=0}^{n}\binom{n%
}{l}x^{l}E_{n-l}^{\left( k\right) }\left( \lambda \right)  \label{eqn6}
\end{equation}%
and%
\begin{equation}
B_{n}^{\left( k\right) }\left( x\mid \lambda \right) =\sum_{l=0}^{n}\binom{n%
}{l}x^{l}B_{n-l}^{\left( k\right) }\left( \lambda \right).  \label{eqaution6}
\end{equation}

From (\ref{eqn1}), (\ref{eqn2}), (\ref{eqn3}), (\ref{equat4}), (\ref{eqn4}),
and (\ref{eqaution7}) we have%
\begin{eqnarray*}
E_{n}^{\left( 1\right) }\left( x\mid \lambda \right) &:&=E_{n}\left( x\mid
\lambda \right) \text{ and }E_{n}^{\left( 1\right) }\left( x\mid 1\right)
:=E_{n}\left( x\mid 1\right) :=E_{n}\left( x\right), \\
B_{n}^{\left( 1\right) }\left( x\mid \lambda \right) &:&=B_{n}\left( x\mid
\lambda \right) \text{ and }B_{n}^{\left( 1\right) }\left( x\mid 1\right)
:=B_{n}\left( x\mid 1\right) :=B_{n}\left( x\right).
\end{eqnarray*}

By (\ref{eqn1}), we easily get%
\begin{equation}
E_{n}^{\left( 0\right) }\left( x\mid \lambda \right) =B_{n}^{\left( 0\right)
}\left( x\mid \lambda \right) =x^{n}.  \label{eqn5}
\end{equation}

Applying derivative operator in the both sides of (\ref{eqaution6}), we have%
\begin{equation}
\frac{d}{dx}B_{n}^{\left( k\right) }\left( x\mid \lambda \right)
=nB_{n-1}^{\left( k\right) }\left( x\mid \lambda \right).  \label{eqn7}
\end{equation}

Using (\ref{eqaution7}), we arrive to 
\begin{equation}
\frac{\lambda B_{n+1}^{(k)}(x+1\mid \lambda )-B_{n+1}^{(k)}(x\mid \lambda )}{%
n+1}=B_{n}^{(k-1)}(x\mid \lambda )\text{,}\quad \text{(see\thinspace \cite%
{Luo1})}\,\text{.}  \label{eqn9}
\end{equation}

The linear operators $\Lambda $ and $D$ on the space of real-valued
differentiable functions are considered as: For $n\in \mathbb{N}$ 
\begin{equation}
\Lambda f(x)=\lambda f(x+1)-f(x)\text{ \textup{and} }Df(x)=\frac{df(x)}{dx}\,%
\text{.}  \label{inteqn100}
\end{equation}

Notice that $\Lambda D=D\Lambda $\thinspace . By (\ref{inteqn100}), we have%
\begin{eqnarray*}
\Lambda ^{2}f\left( x\right) &=&\Lambda \left( \Lambda f\left( x\right)
\right) =\lambda ^{2}f\left( x+2\right) -2\lambda f\left( x+1\right)
+f\left( x\right) \\
&=&\sum_{l=0}^{2}\binom{2}{l}\left( -1\right) ^{l}\lambda ^{l}f\left(
x+l\right) .
\end{eqnarray*}

By continuing this way, we obtain%
\begin{equation*}
\Lambda ^{k}f\left( x\right) =\sum_{l=0}^{k}\left( -1\right) ^{l}\binom{k}{l}%
\lambda ^{l}f\left( x+l\right) .
\end{equation*}

Consequently, we give the following Lemma.

\begin{lemma}
Let $f$ be real valued function and $k\in \mathbb{N}$, we have%
\begin{equation*}
\Lambda ^{k}f\left( x\right) =\sum_{l=0}^{k}\left( -1\right) ^{l}\binom{k}{l}%
\lambda ^{l}f\left( x+l\right) .
\end{equation*}%
In particular,%
\begin{equation}
\Lambda ^{k}f\left( 0\right) =\sum_{l=0}^{k}\left( -1\right) ^{l}\binom{k}{l}%
\lambda ^{l}f\left( l\right) .  \label{eqn8}
\end{equation}
\end{lemma}

Let $P_{n}=\{q(x)\in \mathbb{Q}\lbrack x\rbrack \; |$ $\deg q(x)\leq n\}$ be
the ($n+1$)-dimensional vector space over $\mathbb{Q}$.\thinspace \thinspace
Likely, $\{1,x,\cdots ,x^{n}\}$ is the most natural basis for $P_{n}$.

Additionally, $\{B_{0}^{(k)}\left( x\mid \lambda \right) ,B_{1}^{\left(
k\right) }\left( x\mid \lambda \right) ,\cdots ,B_{n}^{(k)}\left( x\mid
\lambda \right) \}$ is also a good basis for the space $P_{n}$ for our
objective of arithmetical applications of Apostol-Bernoulli polynomials of
higher order.\newline

If $q(x)\in P_{n}$, then $q(x)$ can be written as 
\begin{equation}
q(x)=\sum_{j=0}^{n}b_{j}B_{j}^{\left( k\right) }\left( x\mid \lambda \right) 
\text{.}  \label{equation8}
\end{equation}

Recently, many mathematicians have studied on the applications of \textit{%
polynomials }and $q$-\textit{polynomials} for their finite evaluation
schemes, closure under addition, multiplication, differentiation,
integration, and composition, they are also richly utilized in construction
of their generating functions for finding many identities and formulas (see 
\cite{Araci6}--\cite{Srivastava2}).

In this paper, we discover methods for determining $b_{j}$ from the
expression of $q(x)$ in (\ref{equation8}), and apply those results to
arithmetically and combinatorially interesting identities involving $%
B_{0}^{(k)}\left( x\mid \lambda \right) $, $B_{1}^{(k)}\left( x\mid \lambda
\right) $, $\ldots $, $B_{n}^{(k)}\left( x\mid \lambda \right) $.

\section{Identities on the Apostol-Bernoulli polynomials of higher order}

By (\ref{eqn9}) and (\ref{inteqn100}), we see that 
\begin{equation}
\Lambda B_{n}^{(k)}(x\mid \lambda )=\lambda B_{n}^{(k)}(x+1\mid \lambda
)-B_{n}^{(k)}(x\mid \lambda )=nB_{n-1}^{(k-1)}(x\mid \lambda )\,,
\label{inteqn6}
\end{equation}%
and 
\begin{equation}
DB_{n}^{(k)}(x\mid \lambda )=nB_{n-1}^{(k)}(x\mid \lambda ).  \label{inteqn7}
\end{equation}

Let us assume that $q(x)\in P_{n}$. Then $q(x)$ can be generated by means of 
$B_{0}^{(k)}\left( x\mid \lambda \right)$, $B_{1}^{\left( k\right) }\left(
x\mid \lambda \right)$,\ldots, $B_{n}^{(k)}\left( x\mid \lambda \right)$ as
follows: 
\begin{equation}
q(x)=\sum_{j=0}^{n}b_{j}B_{j}^{(k)}(x\mid \lambda ).  \label{inteqn8}
\end{equation}

Thus, by (\ref{inteqn8}), we get 
\begin{equation*}
\Lambda q\left( x\right) =\sum_{j=0}^{n}b_{j}\Lambda B_{j}^{(k)}(x\mid
\lambda )=\sum_{j=1}^{n}b_{l}jB_{j-1}^{(k-1)}(x\mid \lambda ),
\end{equation*}%
and 
\begin{equation*}
\Lambda ^{2}q(x)=\Lambda \left[ \Lambda q\left( x\right) \right]
=\sum_{j=2}^{n}b_{j}j\left( j-1\right) B_{j-2}^{(k-2)}(x\mid \lambda )\text{.%
}
\end{equation*}

By continuing this way, we have 
\begin{equation}
\Lambda ^{k}q(x)=\sum_{j=k}^{n}b_{j}j\left( j-1\right) \cdots \left(
j-k+1\right) B_{j-k}^{(0)}(x\mid \lambda ).  \label{inteqn9}
\end{equation}

By (\ref{eqn5}) and (\ref{inteqn9}), we see that 
\begin{equation}
D^{s}\Lambda ^{k}q(x)=\sum_{j=k+s}^{n}b_{j}\frac{j!}{\left( j-k-s\right) !}%
x^{j-k-s}  \label{inteqn10}
\end{equation}

Let us take $x=0$ in (\ref{inteqn10}).\thinspace \thinspace\ Then we derive
the following: 
\begin{equation}
\frac{1}{\left( k+s\right) !}D^{s}\Lambda ^{k}q(0)=b_{k+s}.  \label{inteqn11}
\end{equation}

From (\ref{eqn8}) and (\ref{inteqn11}), we have 
\begin{align}
b_{k+s}& =\frac{1}{\left( k+s\right) !}D^{s}\Lambda ^{k}q(0)=\frac{1}{\left(
k+s\right) !}\Lambda ^{k}D^{s}q(0)  \label{inteqn12} \\
& =\frac{1}{\left( k+s\right) !}\sum_{a=0}^{k}\left( -1\right) ^{a}\left( 
\begin{array}{c}
k \\ 
a%
\end{array}%
\right) \lambda ^{a}D^{s}q(a)\text{.}  \notag
\end{align}

Therefore, by (\ref{inteqn8}) and (\ref{inteqn12}), we have the following
theorem.

\begin{theorem}
\label{Theorem 1}For $k\in \mathbb{Z}_{+}$ and $q\left( x\right) \in P_{n}$,
we have%
\begin{equation*}
q(x)=\sum_{j=k}^{n}\left( \frac{1}{j!}\sum_{a=0}^{k}\left( -1\right)
^{a}\left( 
\begin{array}{c}
k \\ 
a%
\end{array}%
\right) \lambda ^{a}D^{j-k}q(a)\right) B_{j}^{(k)}(x\mid \lambda )\text{.}
\end{equation*}
\end{theorem}

Let us take $q(x)=x^{n}\in P_{n}$.\thinspace\ Then we derive that $%
D^{j-k}x^{n}=\frac{n!}{(n-j+k)!}x^{n-j+k}$.\newline

Thus, by Theorem \ref{Theorem 1}, we get 
\begin{equation}
x^{n}=\sum_{j=k}^{n}\left( \frac{1}{j!}\sum_{a=0}^{k}\left( -1\right)
^{a}\left( 
\begin{array}{c}
k \\ 
a%
\end{array}%
\right) \lambda ^{a}\frac{n!}{(n-j+k)!}a^{n-j+k}\right) B_{j}^{(k)}(x\mid
\lambda )\text{.}  \label{inteqn13}
\end{equation}

Therefore, by (\ref{inteqn13}), we arrive at the following corollary.

\begin{corollary}
For $k,n\in \mathbb{Z}_{+}$, we have%
\begin{equation*}
x^{n}=\sum_{j=k}^{n}\left( \frac{1}{j!}\sum_{a=0}^{k}\left( -1\right)
^{a}\left( 
\begin{array}{c}
k \\ 
a%
\end{array}%
\right) \lambda ^{a}\frac{n!}{(n-j+k)!}a^{n-j+k}\right) B_{j}^{(k)}(x\mid
\lambda )\text{.}
\end{equation*}
\end{corollary}


Let $q\left( x\right) =E_{n}^{\left( k\right) }\left( x\right) \in P_{n}$.
Also, it is well known in \cite{Kim6} that%
\begin{equation}
D^{j-k}E_{n}^{\left( k\right) }\left( x\right) =\frac{n!}{\left(
n-j+k\right) !}E_{n-j+k}^{\left( k\right) }\left( x\right) \text{.}
\label{inteqn101}
\end{equation}

By Theorem \ref{Theorem 1} and (\ref{inteqn101}), we get the following
theorem.

\begin{theorem}
For $k,n\in \mathbb{Z}_{+}$, we have%
\begin{equation*}
E_{n}^{\left( k\right) }\left( x\right)
=\sum_{j=k}^{n}\sum_{a=0}^{k}\sum_{l=0}^{n-j+k}\frac{\binom{k}{a}\binom{n-j+k%
}{l}a^{l}\left( -\lambda \right) ^{a}n!}{j!\left( n-j+k\right) !}%
E_{n-j+k-l}^{\left( k\right) }B_{j}^{(k)}(x\mid \lambda )\text{.}
\end{equation*}
\end{theorem}


Let us consider $q\left( x\right) =B_{n}^{\left( k\right) }\left( x\right)
\in P_{n}$. Then we see that%
\begin{equation}
D^{j-k}B_{n}^{\left( k\right) }\left( x\right) =\frac{n!}{\left(
n-j+k\right) !}B_{n-j+k}^{\left( k\right) }\left( x\right) \text{.}
\label{inteqn103}
\end{equation}

Thanks to Theorem \ref{Theorem 1} and (\ref{inteqn103}), we obtain the
following theorem.

\begin{theorem}
For $k,n\in \mathbb{Z}_{+}$, we have%
\begin{equation*}
B_{n}^{\left( k\right) }\left( x\right)
=\sum_{j=k}^{n}\sum_{a=0}^{k}\sum_{l=0}^{n-j+k}\frac{\left( -\lambda \right)
^{a}\binom{k}{a}\binom{n-j+k}{l}a^{l} n!}{j!\left( n-j+k\right) !}%
B_{n-j+k-l}^{\left( k\right) }B_{j}^{(k)}(x\mid \lambda )\text{.}
\end{equation*}
\end{theorem}

Hansen \cite{Hansen} derived the following convolution formula:%
\begin{multline}
\sum_{k=0}^{m}\binom{m}{k}B_{k}\left( x\right) B_{m-k}\left( y\right)= \\
\left( 1-m\right) B_{m}\left( x+y\right) +\left( x+y-1\right) mB_{m-1}\left(
x+y\right) .  \label{eqn10}
\end{multline}

We note that the special case $x=y=0$ of the last identity%
\begin{equation*}
B_{m}=-\frac{\sum_{k=2}^{m-2}\binom{m}{k}B_{k}B_{m-k}}{m+1}
\end{equation*}%
is originally constructed by Euler and Ramanujan (cf. \cite{dilcher}).

Let us now write the following%
\begin{equation}
q\left( x\right) =\sum_{k=0}^{n}\binom{n}{k}B_{k}\left( x\right)
B_{n-k}\left( y\right) \in P_{n}\text{.}  \label{inteqn104}
\end{equation}

By using derivative operator $D^{s}$ in the both sides of (\ref{eqn10}), we
derive%
\begin{gather}
D^{j-k}q\left( x\right) =\left( 1-n\right) \frac{n!}{\left( n-j+k\right) !}%
B_{n-j+k}\left( x+y\right)  \label{inteqn105} \\
+\left( x+y-1\right) \frac{n!}{\left( n-j+k-1\right) !}B_{n-j+k-1}\left(
x+y\right)  \notag \\
+\left( j-k\right) \frac{n!}{\left( n-j+k\right) !}B_{n-j+k}\left( x+y\right)
\notag
\end{gather}

By Theorem \ref{Theorem 1}, (\ref{inteqn104}) and (\ref{inteqn105}), we
arrive at the following theorem.

\begin{theorem}
For $k,n\in \mathbb{Z}_{+}$, we have%
\begin{gather*}
\sum_{k=0}^{n}\binom{n}{k}B_{k}\left( x\right) B_{n-k}\left( y\right) \\
=\sum_{j=k}^{n}\frac{1}{j!}\sum_{a=0}^{k}\left( -1\right) ^{a}\left( 
\begin{array}{c}
k \\ 
a%
\end{array}%
\right) \lambda ^{a}\{\left( 1-n\right) \frac{n!}{\left( n-j+k\right) !}%
B_{n-j+k}\left( a+y\right) \\
+\left( a+y-1\right) \frac{n!}{\left( n-j+k-1\right) !}B_{n-j+k-1}\left(
a+y\right) \\
+\left( j-k\right) \frac{n!}{\left( n-j+k\right) !}B_{n-j+k}\left(
a+y\right) \}B_{j}^{(k)}(x\mid \lambda )\text{.}
\end{gather*}
\end{theorem}

Dilcher \cite{dilcher} introduced the following interesting identity:%
\begin{equation*}
\sum_{k=0}^{n}\binom{n}{k}E_{k}\left( x\right) E_{n-k}\left( y\right)
=2\left( 1-x-y\right) E_{n}\left( x+y\right) +2E_{n+1}\left( x+y\right) 
\text{.}
\end{equation*}

Let $\sum_{k=0}^{n}\binom{n}{k}E_{k}\left( x\right) E_{n-k}\left( y\right)
\in P_{n}$, then we write that 
\begin{equation}
q\left( x\right) =\sum_{k=0}^{n}\binom{n}{k}E_{k}\left( x\right)
E_{n-k}\left( y\right) .  \label{inteqn108}
\end{equation}

By (\ref{inteqn108}), we have%
\begin{gather*}
D^{j-k}q\left( x\right) =2\{\frac{n!}{\left( n-j+k\right) !}\left(
1-x-y\right) E_{n-j+k}\left( x+y\right) \\
-\left( j-k\right) \frac{n!}{\left( n-j+k+1\right) !}E_{n-j+k+1}\left(x+y%
\right) \\
+\frac{\left( n+1\right) !}{\left( n+1-j+k\right) !}E_{n+1-j+k}\left(
x+y\right) \}\text{.}
\end{gather*}

As a result of the last identity and Theorem \ref{Theorem 1}, we derive the
following.

\begin{theorem}
The following equality holds:%
\begin{gather*}
\sum_{k=0}^{n}\binom{n}{k}E_{k}\left( x\right) E_{n-k}\left( y\right) \\
=2\sum_{j=k}^{n}\frac{1}{j!}\sum_{a=0}^{k}\left( -1\right) ^{a}\left( 
\begin{array}{c}
k \\ 
a%
\end{array}%
\right) \lambda ^{a}\{\frac{n!}{\left( n-j+k\right) !}\left( 1-x-y\right)
E_{n-j+k}\left( x+y\right) \\
-\left( j-k\right) \frac{n!}{\left( n-j+k+1\right) !}E_{n-j+k+1}\left(
x+y\right) \\
+\frac{\left( n+1\right) !}{\left( n+1-j+k\right) !}E_{n+1-j+k}\left(
x+y\right) \}B_{j}^{(k)}(x\mid \lambda )\text{.}
\end{gather*}
\end{theorem}

\begin{remark}
Throughout this paper when we take $\lambda =1$, our results can easily be
related to Bernoulli polynomials of higher order.
\end{remark}

\begin{remark}
Theorem \ref{Theorem 1} seems to be plenty large enough for obtaining
interesting identities related to special functions in connection with
Apostol-Bernoulli polynomials of higher order.
\end{remark}

\section*{Acknowledgements}

The authors wish to thank Bernd Kellner for his valuable suggestions for the
present paper.

\end{document}